\documentclass[12pt,twoside,reqno]{amsart}
\usepackage{amsmath}
\usepackage{amsfonts}
\usepackage{amssymb}
\usepackage{color}
\newcommand{\NN}{\mathbb N} 
\newcommand{\RR}{\mathbb R}

\definecolor{dmagenta}{rgb}{0.8,0,0.8}

\newtheorem{theorem}{Theorem}[section]

\newtheorem{proposition}[theorem]{Proposition}
\newtheorem{definition}[theorem]{Definition}
\newtheorem{remark}[theorem]{Remark}
\newcommand\QED{\hfill $\square$}
\begin{document}
\renewcommand{\baselinestretch}{1.18}
\title{Global existence for a hydrogen storage model\\[0.1cm]
with full energy balance}

\author{Elena Bonetti}
\address{Dipartimento di Matematica ``F.~Casorati'', Universit\`a di Pavia, Via Ferrata 1, I-27100~Pavia,~Italy}
\email{elena.bonetti@unipv.it}
\author{Pierluigi Colli}
\address{Dipartimento di Matematica ``F.~Casorati'', Universit\`a di Pavia, Via Ferrata 1, I-27100~Pavia,~Italy}
\email{pierluigi.colli@unipv.it}
\author{Philippe Lauren\c cot}
\address{Institut de Math\'ematiques de Toulouse, CNRS UMR~5219, Universit\'e de Toulouse, F--31062 Toulouse Cedex 9, France} 
\email{laurenco@math.univ-toulouse.fr}
\date{\today}

\begin{abstract}
A thermo-mechanical model describing hydrogen storage by use of metal hydrides has been recently 
proposed in \cite{bfl} 
describing the formation of hydrides using the phase transition approach. 
By virtue of the laws of continuum thermo-mechanics,
the model leads to a phase transition problem in terms of three state variables: the 
temperature, the phase parameter representing the fraction of one solid phase, and the pressure, 
and is derived within a generalization of the principle of 
virtual powers proposed by Fr\'emond \cite{Fremond}, 
accounting for micro-forces, responsible for the phase transition, 
in the whole energy balance of the system. 
Three coupled nonlinear partial differential equations combined with initial and boundary 
conditions have to be solved. The main difficulty in investigating the resulting system of 
partial differential equations relies on the presence of the squared time  
derivative of the order parameter in the energy balance equation, and actually this term was 
neglected in the analysis performed in \cite{bfl}. Here, the global existence 
of a solution to the full problem is proved by exploiting known and sharp estimates on parabolic equations
with right hand side in $L^1.$ Some complementary results on stability and steady state solutions are  
also given.  
\end{abstract}

\maketitle

\renewcommand{\theequation}{\thesection.\arabic{equation}}

%
%
\pagestyle{myheadings}
\newcommand\testopari{\sc{Elena Bonetti --- Pierluigi Colli --- Philippe Lauren\c cot}}
\newcommand\testodispari{\sc{Global existence for a hydrogen storage model}}
\markboth{\testopari}{\testodispari}

\section{Introduction}
\label{intro}
\setcounter{equation}{0}

Let us consider a bounded smooth open set $\Omega \subset \RR^N$, $1\le N\le 3$, with boundary $\Gamma$,
and some final time $T>0$. In the space-time cylindrical domain $Q:=\Omega\times (0,T)$ 
we deal with the following evolution problem
\begin{eqnarray}
&\partial_t e - \Delta \theta = -h(\theta) \partial_t \chi  + \mu\, (\partial_t \chi)^{\, 2} ,&\label{1}\\
&e = \psi(\theta,\chi):= \theta -\chi (h(\theta)- \theta h'(\theta)),&\label{2}\\
&\mu\, \partial_t \chi  - \Delta \chi + \xi = h(\theta)- \log u ,&\label{3}\\
&\xi \in \beta (\chi): = \partial I_{[0,1]} (\chi) + \log (1+\chi),&\label{4}\\
&\partial_t u - \Delta p = 0 ,&\label{5}\\
&p=u(1+\chi),&\label{6}
\end{eqnarray} 
with the boundary conditions 
\begin{equation}
 \partial_n \theta = \partial_n \chi= \partial_n p + \gamma p= 0\quad(\gamma\geq0), \label{7}
\end{equation}
in $\Sigma=\Gamma\times (0,T)$, and the initial conditions 
\begin{equation}
 e (\, \cdot\,, 0) = e_0 , \quad 
 \chi (\, \cdot \,, 0) = \chi_0 , \quad 
 u (\, \cdot\,, 0) = u_0  \quad \hbox{ in } \Omega.
\label{8}
\end{equation}
The system \eqref{1}--\eqref{8} originates from a thermo-mechanical model introduced in \cite{bfl} and
describing hydrogen storage by use of metal hydrides. Hydrogen storage, which turns out a challenging subject in
energetic and industrial applications, basically implies a reduction in the enormous volume of hydrogen gas: so pressure plays a fundamental role in the evolution of the phenomenon. There are basically 
six methods to store hydrogen reversibly with a high volumetric 
and gravimetric density: one of them is by metal hydrides. This technique
exploits the possibility of many metals to absorb hydrogen: such metals and alloys are able to react 
spontaneously with hydrogen and they can store a large amount of it. These materials, either a defined compound 
or a solid solution, are designed as metal hydrides: their use in the hydrogen storage is of interest in terms 
of safety, global yield, and long-time storage. 

Following a usual approach in thermodynamics, Bonetti, Fr\'emond and Lexcellent ~\cite{bfl} derived the governing equilibrium equations (actually simplifying some physical assumptions, but ensuring thermodynamical consistency). The model has next been resumed and extended in \cite{chio}. 

Let us comment about the physical meaning and derivation of equations \eqref{1}-\eqref{5}. As already mentioned, metals are able to absorb hydrogen atoms and combine with them to form solid solutions. The existence of two solid phases is supposed and the occurrence of one phase with respect to the other depends on the pressure $p$ of the hydrogen.  The volume fraction $\chi$ of one of the phases is taken as a state quantity. Then, assuming that no voids appear in the mixture, the volume fraction of the other hydrogen phase is simply given by $1-\chi$. Hence, $\chi$ has to satisfy the relation $0\leq \chi \leq 1,$ which is ensured by the inclusion \eqref{4}  due to the presence of the subdifferential of the indicator function $I_{[0,1]}$ of the interval $[0,1]$ as a component of the maximal monotone graph $\beta $. Recall that the indicator function $I_{[0,1]}$ is given by $I_{[0,1]}(\chi)=0$ if $\chi\in[0,1],$ and  $I_{[0,1]}(\chi)=+\infty$ otherwise. Note that $\beta$ also contains the contribution of the smooth (when restricted to $[0,1]$) increasing function $\log (1+\chi)$ (in our notation, $\log r = s$ if $r= e^s$). The graph  $\beta $ is actually multivalued and maximal monotone (see, e.g., \cite{Brezis, Barbu}), with $ \xi \in \beta(\chi)$ if and only if
$$
\xi - \log (1+\chi) \ \left\{
\begin{array}{ll}
\displaystyle
\leq \, 0 \   &\hbox{if } \ \chi=0  
\\[0.1cm]
= \, 0 \   &\hbox{if } \ 0 <\chi <1  
\\[0.1cm]
\geq \, 0 \  &\hbox{if } \  \chi=1  
\\[0.1cm]
\end{array}
\right. .
$$
Besides the phase parameter $\chi$ along with its gradient $\nabla\chi$ accounting for local interactions between the different phases, the other state variables of the model are the absolute temperature $\theta$ and the hydrogen pressure $p$. Constitutive relations for the state quantities are chosen in such a way that the principles of thermodynamics are satisfied. 

The first equation \eqref{1} is derived from the first principle of thermodynamics. It is obtained generalizing the principle of virtual powers as proposed by Fr\'emond in \cite{Fremond}, i.e. including internal power of micro-forces (which are responsible for phase transitions) in the energy balance. In particular, note that in the right hand side, we have the quadratic nonlinearity $(\partial_t \chi )^2$ standing for the internal power of dissipative micro-forces. The auxiliary variable $e$ is prescribed  in terms of $\theta$ and $\chi$  by relation \eqref{2} via the smooth function $\psi$. In fact, $\psi$ satisfies all the good properties we need (in particular, invertibility with respect to the first variable) since the function $h$ is smooth and small enough (see \cite{bfl, chio}): more precisely, we assume that there exists a constant  $c_h > 1$ such that  
\begin{eqnarray}
& \displaystyle h  \in C^2(\RR) \cap W^{2,\infty} (\RR), \quad \ \| h \|_{W^{2,\infty} (\RR)} + \sup_{r\in\RR} |r h'(r)| \leq c_h &\label{12}\\
& \displaystyle\frac1{c_h} \leq 1+r h''(r)s \leq c_h \quad \hbox{for all } \ r\in\RR \hbox{ and } s\in [0,1]. &\label{13}
\end{eqnarray}

Equation \eqref{3}, describing the evolution of the phase parameter, comes 
from the already mentioned generalization of the principle of virtual power in which micro-movements 
(related to the evolution of $\chi$ as they are linked to changes of micro-structure) 
are included \cite{Fremond}. Equation \eqref{5} is deduced from the usual continuity 
equation, in which the velocity is assumed to be proportional to the gradient of the pressure. The structure of equation \eqref{3} (along with \eqref{4}) is by now quite standard in the 
framework of phase field models: observe however that in the right hand side one has to 
handle the singular term $\log u$ (or alternatively $\log p$). Then, one has to recover 
admissibility and enough regularity for $\log u$ from both \eqref{5}--\eqref{6} and the 
initial datum $u_0$, that has to be strictly positive. The same can be repeated for $e_0$ 
if we want to guarantee that the absolute temperature $\theta $ remains positive during 
the evolution, as it should be: but concerning \eqref{1}, we can be more flexible and 
develop our existence theory also in the general case, without sign restriction.

About boundary conditions \eqref{7}, we explain that $\partial_n $ denotes the 
outward normal derivative on the boundary. Hence, Neumann homogeneous boundary conditions are set for 
$\theta $ and $\chi$,  while the condition for $p$ states that the hydrogen flux through the boundary is somehow 
proportional to the value of the interior pressure on $\Gamma$. Anyhow, in our contribution 
we treat the case $\gamma =0 $ too, corresponding to the case when the system is perfectly isolated. 

Now, let us briefly discuss the presence of dissipative quadratic 
nonlinearities and the coupling of phase transitions with pressure. Concerning phase 
changes with microscopic motions, there is a comprehensive literature on a class of models, 
some of them originating from the theory developed by Fr\'emond \cite{Fremond}. In this respect, we first refer the reader to \cite{BFL}: the resulting system of phase field type is characterized by the 
presence of new nonlinearities, like  $(\partial_t \chi )^2$ and $\theta 
\partial_t\chi$, which were not present in the classical formulation of phase field 
systems (not accounting for microscopic stresses). Several authors have dealt with 
this kind of problems and various situations have been analyzed. However, mainly for 
analytical difficulties due to the presence of nonlinearities,  to our knowledge there 
is no global in time well-posedness result for the complete related system 
in the 3D (or 2D) case. A global existence result is proved in the 1D setting \cite{LSS,LS} or for a non-diffusive phase evolution \cite{CLSS}. Other results have been obtained for some regularized versions of the problem \cite{B1}. Weak solutions in terms of energy conservation (of Feireisl type) have been recently constructed in \cite{FPR}. In addition, if we focus on the presence of the pressure, and 
thus on the continuity equation, note that some recent papers have dealt with phase 
transitions allowing different densities of the phases or adding the possibility of 
voids (see, for instance, \cite{BFF,RF}, where internal constraints on 
the pressure are given through indicator functions).

In our system, the unpleasant nonlinearity in the right hand side of \eqref{1} includes $\mu\, (\partial_t \chi )^{\, 2}$ with $\mu $ being a positive parameter. Note that the same $\mu $ appears as the 
coefficient of $\partial_t \chi  $ in  \eqref{3}. Since our aim is showing the existence of a possibly weak solution to 
the initial boundary value problem \eqref{1}--\eqref{8}, the presence of $\mu\, (\partial_t \chi )^{\, 2}$ in 
\eqref{1} creates, as we have already pointed out,  a difficulty from the mathematical point of view.  For this reason, the
authors of \cite{bfl} suppressed it in the first analysis they carried out after the derivation of the model. That term is instead maintained in \cite{chio}, where equation \eqref{3} is replaced by 
\begin{equation}
\mu \partial_t \chi  - \nu\Delta \partial_t \chi - \Delta \chi + \xi = h(\theta)- \log u , \label{3reg}
\end{equation}
where $\nu$ stands for a small positive coefficient. The additional term $- \nu\Delta \partial_t \chi $ actually 
contributes for a better regularity of $\partial_t \chi $, and consequently the nonlinearity $\mu\, (\partial_t \chi )^{\,2}$ to be handled in \eqref{1} turns out to be smoother. 

In the present contribution we can overcome the mentioned difficulty by exploiting an argument by 
Boccardo and Gallou\" et \cite{bg} for the study of second-order nonlinear parabolic equations with 
right hand side in $L^1$.  We prove the global-in-time existence of a weak solution to \eqref{1}--\eqref{8}
by taking (in Section~\ref{appr}) an approximation of the problem  based on \eqref{3reg}, 
with a small parameter $\nu$ intended to go to zero in the limit procedure. 
Then we derive a priori estimates in Section~\ref{apriori} and finally perform 
the passage to the limit: this is the subject of Section~\ref{limit} which concludes the existence proof. 
Section~\ref{complement} is devoted to the proof of some complementary results stated in Section~\ref{mainres}
(along with the main existence theorem) and regarding positivity of the variable $\theta$, 
stability estimates in the case $\gamma>0$ and stationary solutions for the case $\gamma =0$. 

\section{Main result}
\label{mainres}
\setcounter{equation}{0}

In this section, we make precise assumptions on the data, 
specify our notion of solution to problem \eqref{1}--\eqref{8} 
and state our existence result.
For convenience we set
\begin{equation}
  V := H^1(\Omega)
  \quad H := L^2(\Omega) , \quad
  W := 
\left\{v\in H^2(\Omega):\ \partial_n v = 0 
\;\,{\textrm{on}\,\; \Gamma} \right\} ,
  \label{defspazi}
\end{equation}
and endow these spaces with their standard norms,
for which we use {the self-explanato\-ry} notations  $\|\,\cdot\,\|_H$, $\|\,\cdot\,\|_V$, and
$\|\,\cdot\,\|_W$ {(the norm in $L^2(\Omega;\RR^N)$ being also denoted by $\|\,\cdot\,\|_H$)}.
We remark that the embeddings $W\subset V\subset H$ are compact,
because $\Omega$ is bounded and smooth.
Since $V$ is dense in~$H$, we can identify $H$  with a subspace of $V'$ in the usual~way
(i.e.,~{so as to have} that $\langle u,v \rangle=(u,v)_H$
for every $u\in H$ and $v\in V$); the embedding $H\subset V'$ is also compact. Let $T>0$ and set $Q= 
\Omega\times (0,T)$ and $\Sigma=\Gamma\times (0,T)$. 

We next recall that
\begin{equation}
 \mu >0 , \quad \gamma \geq 0, \quad h:\RR\to\RR \ \hbox{ satisfies \eqref{12}--\eqref{13}}
  \label{hyp1}
\end{equation}
and introduce the operators $A_\gamma, \, A :V\to V'$ specified by
\begin{equation}
\langle A_\gamma v , w \rangle := \int_\Omega \nabla v \cdot \nabla w \, dx + 
\gamma \int_\Gamma v\, w \, dx
\quad \forall\ v,w \in V, \quad \ A\equiv A_0 \, .
\label{15}
\end{equation}
In view of the definition of $\psi :\RR \times [0,1] \to \RR$ given in \eqref{2}, it is 
straightforward to check that the partial derivatives $\partial_1 \psi, \, \partial_2 \psi $ 
(with respect to first and second variable, respectively) satisfy 
\begin{equation}
\begin{split}
\displaystyle\frac1{c_h} \leq \partial_1 \psi (r,s) = 1+ r h''(r)s  \leq c_h , \quad 
|\partial_2 \psi (r,s) | & \leq |h(r)| + |r h'(r)|  \leq c_h \qquad \\
& \hbox{for all } \ r\in\RR  \hbox{ and } s\in [0,1]. 
\end{split}
\label{psi1}
\end{equation}
In particular, by looking at \eqref{2} again, we remark that for all $e\in \RR$ and $\chi\in [0,1]$ 
there exists a unique $\theta \in \RR$ fulfilling $e=\psi (\theta, \chi)$, and we denote by $\psi^{-1}$ 
the inverse function of $\psi$ with respect to the first variable. Clearly, the partial derivatives of $\psi^{-1}$ satisfy 
\begin{equation}
\begin{split}
\displaystyle\frac1{c_h} \leq \partial_1 \psi^{-1} (r,s) \leq c_h , \quad 
|\partial_2 \psi^{-1} (r,s) | \leq &  |\partial_1\psi^{-1} (r,s)\partial_2\psi (\psi^{-1} (r,s), s)|  \leq 
(c_h)^{2} \qquad \\
& \hbox{for all } \ r\in\RR  \hbox{ and } s\in [0,1]. 
\end{split}
\label{psi2}
\end{equation}
In order to develop our existence theory, due to the presence of the term $\mu\, (\partial_t \chi )^{\,2}$ in \eqref{1} 
we should give a meaning to \eqref{1} also when the right hand side is just in $L^1(Q)$. To this aim, 
let us extend the operator $A$ in \eqref{15} from $W^{1,q}(\Omega)$ to $(W^{1,q'}(\Omega))'$, 
\begin{equation}
\begin{split}
\langle A v , w \rangle := \int_\Omega \nabla v \cdot \nabla w \, dx
& \quad \mbox{ for all }\;\; (v,w) \in W^{1,q}(\Omega) \times W^{1,q'}(\Omega), \hskip2cm\\
& \hbox{where }  1<q < \frac{N+2}{N+1}, \quad \frac{1}{q} +\frac{1}{q'}=1,  \quad q' > N+2.
\end{split}
\label{Aext}
\end{equation}

We now introduce the proper convex lower semicontinuous function $\hat\beta :\RR \to [0,+\infty )$ defined by 
\begin{equation}
\label{defbetahat}
\hat\beta (r) =  \left\{
\begin{array}{ll}
\displaystyle
\int_0^r \log(1+s)ds  \   &\hbox{if } \ 0 \leq r \leq 1  
\\[0.4cm]
+\infty \  &\hbox{elsewhere }  
\\
\end{array}
\right. ,
\end{equation}
so that  $\hat\beta(r)\ge \hat\beta(0)=0$ for $r\in\RR$, and its subdifferential is nothing but 
$\beta$ (cf. \eqref{4}). Hence, let us consider  the functionals (cf. \cite[formulas (58)--(66)]{chio})
\begin{equation}
\label{defJ}
J_H (v) = \left\{
\begin{array}{ll}
\displaystyle
\int_\Omega \hat\beta (v) \,dx  \   &\hbox{if } \ v\in H \ \hbox{ and } \hat\beta (v)\in L^1(\Omega)    
\\[0.4cm]
+\infty \  &\hbox{if } \ v\in H \ \hbox{ but } \hat\beta (v)\not\in L^1(\Omega)
\\
\end{array}
\right. , \ \quad \ J(v) = J_H (v) \quad \hbox{ if } \ v\in V, 
\end{equation}
It is known that the inclusion $w\in \beta (v)$ a.e. in $\Omega $ can be rewritten as $w\in \partial_H J_H (v)$ provided $(v,w)\in H\times H$ or reinterpreted in terms of the subdifferential $\partial_{V,V'} J : V \to 2^{V'}$ (different from the previous subdifferential $\partial_{H} J_H : H \to 2^{H}$) as $w\in \partial_{V,V'} J (v)$ if $(v,w)\in V\times V'$. Thus, this last inclusion can be seen as an extension of the previous ones whenever $v\in V$. We note in particular that
\begin{equation}
\partial_H J_H (v) \equiv H\cap \partial_{V,V'} J (v) \quad\mbox{ for all }\;\;\ v\in V,
\label{sotdif}
\end{equation}
and announce that for our convenience we will use $\partial_{V,V'} J $ in the variational formulation of \eqref{1}--\eqref{8} to account for relation \eqref{4}.

Concerning equation \eqref{5}, we point out that we shall take advantage of a time-integrated version of it
which is obtained with the help of the corresponding initial condition in \eqref{8}, that is,
\begin{equation}
u -\Delta (1*p) = u_0, 
\label{5bis}
\end{equation}
where $ \displaystyle (1*p) (\, \cdot\,, t) = \int_0^t  p (\, \cdot\,, s) ds  \hbox{ for } t\in (0,T).$    

At this point, we can prescribe the conditions on the initial data. We assume that
\begin{equation}
\begin{split}
& e_0 \in L^1 (\Omega), \quad \chi_0 \in V \, \hbox{ with } \,\hat\beta (\chi_0)\in L^1(\Omega), \hskip4cm \\ 
& u_0 \in  L^1(\Omega)\cap V' \,\hbox{ such that } \, \log u_0 \in L^1(\Omega).
\end{split}
\label{35}
\end{equation}
Let us just notice that the requirement on $\chi_0$ exactly says that $\chi_0 \in D(J)$, i.e. $\chi_0$
lies in the domain of $J$, and this obviously entails $0\leq \chi_0 \leq 1 $ a.e. in $\Omega $. The 
assumption on $u_0$ implies that both $u_0$ and $ \log u_0 $ are in $L^1(\Omega)$, whence in particular $u_0$ is strictly positive a.e. in $\Omega $.

Here is our notion of weak solution to the problem \eqref{1}--\eqref{8}.

\begin{definition}
\label{solution}
We say that the sextuple $(e,\theta,\chi,\xi,u,p)$ is a weak solution to problem \eqref{1}--\eqref{8} if 
\begin{eqnarray}
& \displaystyle e , \,\theta \in L^q(0,T;W^{1,q}(\Omega)) \quad \hbox{for some $q$ fulfilling } \ 1<q < \frac{N+2}{N+1},&
\label{29}\\
&\displaystyle  \partial_t e \in L^1(0,T;(W^{1,q'}(\Omega))') \quad \hbox{ with } \ q'=\frac{q}{q-1} > N + 2,
& \label{30}\\
& \chi \in H^1(0,T;H)\cap L^2(0,T;V), & \label{31}\\
& \xi \in L^2(0,T;V') ,& \label{32}\\
& u\in L^\infty(0,T;V') \cap L^2(0,T;H) \ \hbox{ and } \ u >0 \ \hbox{ a.e. in } Q ,  & \label{33}\\
& p \in  L^2(0,T;H)   \ \hbox{ with } \  1*p \in L^\infty(0,T;V)  & \label{34}
\end{eqnarray}
and the following equations and conditions are satisfied:
\begin{eqnarray}
&\partial_t e + A \theta = -h(\theta) \partial_t \chi  + \mu \partial_t \chi ^{\, 2} 
\quad\hbox{in } (W^{1,q'}(\Omega))', \hbox{ a.e. in } (0,T),& \label{22}\\
&e = \psi(\theta,\chi)
\quad\hbox{a.e. in } Q,& \label{23}\\
&\mu \partial_t \chi  + A \chi + \xi = h(\theta) - \log u 
\quad\hbox{in } V', \hbox{ a.e. in } (0,T),& \label{24}\\
&\xi (t) \in \partial_{V,V'} J (\chi(t)) 
\quad\hbox{for a.e. } t\in (0,T),& \label{25}\\
&u + A_\gamma (1*p) = u_0 \quad\hbox{in } V', \hbox{ a.e. in } (0,T),& \label{26}\\
&p=u(1+\chi)\quad\hbox{a.e. in } Q,& \label{27}\\
&e (0) = e_0 \quad\hbox{in } (W^{1,q'}(\Omega))', \ \quad \chi (0) = \chi_0 \quad\hbox{in } H. & \label{28}
\end{eqnarray}
\end{definition}

As a remark, note that the right hand side of \eqref{22} is only in $L^1(Q)$, and we have $L^1(\Omega) \subset (W^{1,q'}(\Omega))'$. Indeed, seeing that $q'>N$, then $W^{1,q'}(\Omega) $ is densely and compactly embedded in $L^\infty (\Omega)$. We also point out that the inclusion \eqref{25} yields
\begin{equation}
0\leq \chi\leq 1 \quad\hbox{a.e. in } Q, \label{25bis}\\
\end{equation}
as a consequence, since $\chi(t) $ has to belong to $D(\partial_{V,V'} J) \subseteq D( J)$ 
for a.e. $t\in (0,T).$ The positivity requirement on $u\in L^2(Q)$ (cf.~\eqref{33}) ensures that $\log u $
is well defined and  measurable as a function: in addition, by comparing terms in \eqref{24} we read that
$\log u \in L^2(0,T;V')$ at least.
The initial conditions  in \eqref{28} make sense due to \eqref{29}--\eqref{30}, which entail 
the continuity of $e$ from $[0,T]$ to  $(W^{1,q'}(\Omega))'$, and \eqref{31}.

It is time for us to state our existence result.
\begin{theorem}
\label{existence}
Under the assumptions \eqref{hyp1} and \eqref{35}, 
there exists a weak solution $(e,\theta,\chi,\xi,u,p)$ to problem \eqref{1}--\eqref{8}  in the 
sense of Definition~\ref{solution}. Moreover, the further regularity properties
\begin{eqnarray}
&  e , \,\theta \in L^\infty(0,T;L^1(\Omega))\cap L^r(Q)\cap L^q(0,T;W^{1,q}(\Omega)) \hskip3cm &
\nonumber \\
&\hskip3cm \displaystyle\hbox{for all $r,\,q $ fulfilling } \ 
1\leq r <\frac{N+2}{N}, \  1\leq q < \frac{N+2}{N+1},&
\label{29bis}\\[0.2cm]
& \chi \in C^0([0,T];V)\cap L^2(0,T;W), & \label{31bis}\\[0.2cm]
& \xi \in L^2(0,T;H) ,& \label{32bis}\\[0.2cm]
& u\in L^\infty(0,T;L^1(\Omega)) \ \hbox{ and } \ \log u \in  L^\infty(0,T;L^1(\Omega))\cap L^2(0,T;V)   & \label{33bis}
\end{eqnarray}
hold true for $(e,\theta,\chi,\xi,u,p)$.
\end{theorem}

\begin{remark}
\label{rem-dopo-def}
\rm In view of \eqref{sotdif} and \eqref{25}, it turns out that \eqref{32bis} ensures the validity of $\xi (t) 
\in \partial_{H} J_H (\chi(t))$ for a.e. $t\in (0,T),$ which in turn yields \eqref{4}, that is 
$\xi \in \beta  (\chi) $ a.e. in $Q$. Thanks to \eqref{31bis} and \eqref{33bis}, we infer from \eqref{24}  that \eqref{3} is solved almost everywhere as well, which is not the case for the other equations \eqref{1} and \eqref{5}. In fact, for that we would need a additional regularity for $(e,\theta)$ and $(u,p)$, respectively. Note however that equation \eqref{5bis} is satisfied a.e. in $Q$ if $u_0 \in H$, since in this case 
\begin{equation}
1*p\in  L^2(0,T;H^2(\Omega));  
\label{1convp}
\end{equation}
indeed, \eqref{1convp} is a consequence of  \eqref{34}, $ \Delta (1*p) \in  L^2(0,T;H)$ 
and known elliptic regularity estimates which exploit the smoothness 
of the domain and of the boundary condition. Anyway, we would like to emphasize that 
\eqref{5} holds true in the sense of distributions in $Q$. 
\end{remark}

The next statement regards positivity of the variable $\theta$, which is consistent with the fact that it represents the absolute temperature.
\begin{theorem}
\label{posit-temp}
Assume that
\begin{equation}
\log \theta_0  \in L^1(\Omega), \quad \hbox{where} \quad \theta_0:= \psi^{-1} (e_0, \chi_0),
\label{hyptheta0}
\end{equation}
in addition to \eqref{hyp1} and \eqref{35}. Then the solution $(e,\theta,\chi,\xi,u,p)$
specified by Theorem~\ref{existence} satisfies 
\begin{equation}
\log \theta  \in L^\infty(0,T;L^1(\Omega))\cap L^2(0,T;V), \quad \frac{\partial_t \chi }{\sqrt\theta} \in L^2(Q) 
\label{reg1}
\end{equation}
and in particular $\theta >0 $ a.e. in $Q$.
\end{theorem}
We also have two results concerning the large time behaviour in  the case $\gamma >0$ and 
the steady state solutions for the problem with $\gamma=0$.
\begin{theorem}
\label{decay}
Let $\gamma>0$ and the assumptions of Theorem~\ref{existence} hold. Then, we have that 
$$ u\in L^2(0,+\infty; L^2(\Omega)) \quad  \hbox{and } \ \|u(t)\|_{V'} \ \hbox{ decays to zero exponentially
as } \ t\to +\infty. $$
\end{theorem}

In connection with this statement, note that $u(t_n + \, \cdot)  \to 0$ 
strongly in $L^2(0,T;H)$ for all $T>0$ and for any sequence $t_n \to +\infty$. In some improper 
way, this kind of convergence to zero of $u$ would entail that the right hand side of \eqref{24} 
tends somehow to  $+\infty$, due to the presence of $- \log u$. 
Thus, we expect that $\chi(t)$ goes to $1$, or becomes $1$ from a certain time, which corresponds to  
one of the two solid hydrogen phases that becomes dominant at a long range. 

\begin{theorem}
\label{stationary}
Let $\gamma= 0$. If $(e,\theta,\chi,\xi,u,p)$ is a stationary weak solution to \eqref{1}--\eqref{8}, then the elements of the sextuple are all constant functions satisfying 
\begin{eqnarray*}
&\hbox{if } \ h(\theta) - \log p >0  \quad &\hbox{then } \ \chi = 1 \ \hbox{ and }\ p = 2 u, \\
&\hbox{if } \ h(\theta) - \log p =0  \quad &\hbox{then } \ 0 \leq \chi \leq  1 \ \hbox{ and }\ p = (1+\chi) u , \\
&\hbox{if } \ h(\theta) - \log p <0  \quad &\hbox{then } \ \chi =  0 \ \hbox{ and }\ p = u.
\end{eqnarray*}
\end{theorem}

\section{Approximation of the problem}
\label{appr}
\setcounter{equation}{0}

In order to prove Theorem~\ref{existence}, we introduce a regularizing term in \eqref{24},
take advantage of the results shown by Chiodaroli in \cite{chio}, prove uniform estimates 
and then pass to the limit with respect to the approximation parameter. 

We also regularize initial data and deal with a sequence of strictly 
positive values $\{ \gamma_n \}_{n\in \NN}$ to be used in the boundary condition for $p$ (cf.~\eqref{7}): more precisely, we assume that
\begin{equation}
\gamma_n>0\,, \quad \gamma_n\to \gamma,  \ \hbox{ whence } \ A_{\gamma_n} \to A_\gamma \ \hbox{ in } \ \mathcal{L}(V; V'), \ \hbox{ as } \ 
n\to \infty.
\label{convgamma}
\end{equation}
Besides, we require the following conditions for the approximating initial data:
\begin{eqnarray}
&\chi_{0n} \in W \cap D(\partial_{V,V'} J ),\ \hbox{ and consequently } \  
\,\hat\beta (\chi_{0n})\in L^1(\Omega),  &
\label{37}\\ 
&e_{0n} ,\, \theta_{0n} \in V, \ \hbox{ where } \ \theta_{0n} := \psi^{-1} (e_{0n}, \chi_{0n}), &
\label{36}\\
&\displaystyle u_{0n} ,\, p_{0n} \in V,
\ \hbox{ where } \ p_{0n} :=  u_{0n} (1+\chi_{0n}), \
\hbox{ and } \, \log u_{0n} \in L^1(\Omega), \ \, \frac{1}{p_{0n}} \in H. &
\label{38}
\end{eqnarray}
Clearly, $e_{0n} ,\, \chi_{0n}, \, u_{0n} $ are smoother than the data $e_{0} ,\, \chi_{0}, \, u_{0} $ 
satisfying assumption \eqref{35}. Morever, we also ask that
\begin{eqnarray}
&e_{0n} \to  e_0 \ \hbox{ in } \ L^1(\Omega), \quad \chi_{0n} \to  \chi_0 \ \hbox{ in } \ V, \quad 
u_{0n} \to  u_0 \ \hbox{ in } \ V' \quad \hbox{ as } \ 
n\to \infty, & \label{39}\\
&\displaystyle \| \chi_{0n} \|_W^2 \leq C\,  n, \quad \| u_{0n}- \log u_{0n} \|_{L^1(\Omega)} \leq C \quad \hbox{ for all }
n\in\NN.&
\label{39bis}
\end{eqnarray}
Here and in the sequel, $C $ denotes a positive constant, possibly depending on data but independent of $n$, 
that may change from line to line and even in the same chain of inequalities. Let us also note that we do not 
need to prescribe any boundedness property for $\{ \hat\beta (\chi_{0n})\} $ in $L^1(\Omega)$: 
indeed, this is ensured by \eqref{37}, which implies $0\leq \chi_{0n} \leq 1 $ a.e. in $\Omega $ and consequently 
(cf.~\eqref{defbetahat})
$$
J(\chi_{0n})= \big\| \hat\beta (\chi_{0n}) \big\|_{L^1(\Omega)}  \leq \int_\Omega |\chi_{0n}|  
\, \log 2 \, dx \leq |\Omega|.
$$
Moreover, the strong convergence of $\chi_{0n} $ to $\chi_0$ stated in \eqref{39} also entails 
\begin{equation}
\label{39ter} 
J(\chi_{0n}) \to  J(\chi_0 )\quad \hbox{ as } \ n\to \infty. 
\end{equation}

Let us now briefly sketch how to construct the sequence of approximating initial conditions $(e_{0n} ,\, \chi_{0n}, \, u_{0n}) $ satisfying \eqref{37}--\eqref{39ter}. First, there is no problem in constructing  a sequence $\{ e_{0n} \}$ that fulfills \eqref{36} and converges to $e_0$ in $L^1(\Omega)$: it suffices to extend $e_0$ outside of $\Omega$ by zero and then regularize it by taking a convolution product over $\RR^N$ with a sequence of mollifiers. 

As  for $\{\chi_{0n}\}$, we consider the sequence of solutions to
\begin{equation}
\label{chizeron}
\chi_{0n} + \frac1n A \chi_{0n} = \chi_0 \quad \hbox{in } \, V' \,. 
\end{equation}
Since $\chi_0\in V$, we have $\chi_{0n} \in W $  and $0\leq \chi_{0n} \leq 1 $ a.e. in $ \Omega $ by the maximum
principle. Moreover, we have that $\chi_{0n} \in D(\partial_{V,V'} J )$: indeed, 
 $\log (1+\chi_{0n} ) \in \partial_{V,V'} J (\chi_{0n}) $ as 
\begin{eqnarray*}
J(\chi_{0n})= \int_\Omega \hat\beta (\chi_{0n})\, dx \leq 
\int_\Omega \log (1+\chi_{0n} )\,(\chi_{0n}- \chi_0)\, dx +
\int_\Omega \hat\beta (\chi_{0})\, dx \\
= \langle \log (1+\chi_{0n} ) , \chi_{0n}-\chi_0 \rangle + J(\chi_0).
\end{eqnarray*}
By the choice \eqref{chizeron}, the first inequality in \eqref{39bis} holds true. To check it, it 
suffices to test \eqref{chizeron} by $\chi_{0n}$ first and then by $A\chi_{0n} $,  
integrating by parts in the right hand side.
Furthermore, in singular perturbations like \eqref{chizeron} 
a standard argument allows to deduce the strong convergence 
$\chi_{0n} \to \chi_0$ by proving weak convergence plus convergence of norms. 

The same kind of approximation can be assumed for $u_{0}$. In fact, let us take the solution $u_{0n}$ to 
\begin{equation}
\label{uzeron}
u_{0n} + \frac1n A u_{0n} = u_0 +\frac1n \quad \hbox{in } V',
\end{equation}
with the additional $1/n$ in the right hand side. It is straightforward to see that 
$u_{0n} \in V $ and $u_{0n} \geq 1/n $ a.e. in $\Omega$, the latter following from the comparison principle. Then \eqref{38} follows 
and $1/p_{0n} $ is even in $V\cap L^\infty(\Omega)$.
Further, if we test \eqref{uzeron} by $(1 -1/u_{0n})$ and observe that
$r\mapsto (1-1/r)$ is the derivative of the positive and convex function $r\mapsto 
r-\log r$, we find out that 
\begin{eqnarray*}
&&\int_\Omega (u_{0n} - \log u_{0n})\,dx - \int_\Omega (u_{0} - \log u_{0}) \, dx
\leq \int_\Omega \left(1-\frac1{u_{0n}}\right)(u_{0n} - u_{0})\, dx\\
&&\hskip1cm \leq \int_\Omega \left(1-\frac1{u_{0n}}\right)(u_{0n} - u_{0}) \, dx +\frac1n 
\int_\Omega |\nabla \log u_{0n}|^2 \, dx \leq 
\int_\Omega \frac1n\left(1-\frac1{u_{0n}}\right)\, dx
 \leq \frac{|\Omega|}n
\end{eqnarray*}
as well as
$$
0 \leq \int_\Omega (u_{0n} - \log u_{0n}) \, dx 
\leq  \int_\Omega (u_{0} - \log u_{0})\, dx +|\Omega|,
$$
and the property $\| u_{0n}- \log u_{0n} \|_{L^1(\Omega)} \leq C$ turns out from \eqref{35}.

The following proposition is a consequence of \cite[Theorems~2.3 and~2.5]{chio} and it holds for all $n\geq 1$.
\begin{proposition}
\label{app-sol}
Under the assumptions \eqref{hyp1} and \eqref{37}--\eqref{38},  there exists
a sextuple $$(e_n,\theta_n,\chi_n,\xi_n,u_n,p_n)$$ fulfilling 
\begin{eqnarray}
& \displaystyle e_n , \, \theta_n \in  H^1 (0,T;H)\cap L^\infty (0,T; V)   ,&
\label{44}\\
& \chi_n \in W^{1,\infty}(0,T;V)\cap L^\infty(0,T;W)\cap L^\infty (Q), & \label{45}\\
& \xi_n \in L^2(0,T;V') , & \label{46}\\
& u_n\in H^1(0,T;H) \cap L^\infty(0,T;V) \ \hbox{ and } \ \log u_n \in  L^\infty(0,T;L^1(\Omega))\cap L^2(0,T;V), & \label{47}\\
& p_n \in  H^1(0,T;H) \cap L^\infty(0,T;V) \cap L^2(0,T;H^2(\Omega)) & \label{48}
\end{eqnarray}
and the following equations and conditions are satisfied:
\begin{eqnarray}
&\partial_t e_n + A \theta_n = -h(\theta_n) \partial_t  \chi_n + \mu\, (\partial_t \chi_n)^2 
\quad\hbox{in } V', \hbox{ a.e. in } (0,T),& \label{22n}\\
&e_n = \psi(\theta_n,\chi_n) \quad\hbox{a.e. in } Q,& \label{23n}\\
&\displaystyle \mu\, \partial_t\chi_n + \frac1n A (\partial_t \chi_n) + A \chi_n + \xi_n = 
h(\theta_n) - \log u_n \quad\hbox{in } V', \hbox{ a.e. in } (0,T),& \label{24n}\\
&\xi_n (t) \in \partial_{V,V'} J (\chi_n(t)) 
\quad\hbox{for a.e. } t\in (0,T),& \label{25n}\\
&u_n + A_{\gamma_n} (1*p_n) = u_{0n} \quad\hbox{in } V', \hbox{ a.e. in } (0,T),& \label{26n}\\
&p_n=u_n(1+\chi_n)\quad\hbox{a.e. in } Q,& \label{27n}\\
&e_n (0) = e_{0n} \quad\hbox{in } V', \ \quad \chi_n (0) = \chi_{0n} \quad\hbox{in } H. & \label{28n}
\end{eqnarray}
\end{proposition}

After the statement, we note that the regularity properties \eqref{44}--\eqref{48} entail 
\eqref{29}--\eqref{34} and \eqref{29bis}--\eqref{31bis}, \eqref{33bis}. 
On the other hand, \eqref{32bis} is not granted by \eqref{46}: in fact, from \eqref{24n} it just follows that 
$$n^{-1} A (\partial_t \chi_n) + \xi_n  \in L^\infty(0,T;H) +L^2(0,T;V)$$
(cf.~\eqref{45} and \eqref{47}), but  $A (\partial_t \chi_n)$ is only in $L^\infty(0,T;V')$.
An additional remark concerns equation \eqref{26n} which can be  rewritten for  $u_n$ and  $p_n$ as 
\begin{equation}
\partial_t u_n + A_{\gamma_n} p_n = 0 \quad\hbox{in } V', \hbox{ a.e. in } (0,T) \label{50}
\end{equation}
along with the initial condition 
\begin{equation}
u_n (0) = u_{0n} \quad\hbox{in } H. \label{51}
\end{equation}
Note that \eqref{50} and \eqref{47}--\eqref{48} imply that $u_n$ and $p_n$ solve \eqref{5} and the 
related boundary condition almost everywhere. The same can be established for \eqref{22n} which yields the analog of
\eqref{1} a.e. in $Q$: indeed, \eqref{44}--\eqref{45} entail $A\theta_n \in L^2(0,T;H)$ so that
$\theta_n\in L^2(0,T;W)$ too. As a last remark, we point out that \eqref{44}, \eqref{23n} and \eqref{12} allow 
us to rewrite \eqref{22n} as 
\begin{equation}
\left( 1 +  \chi_n \theta_n h''(\theta_n) \right) \partial_t \theta_n 
+ A \theta_n = -\theta_n h'(\theta_n) \partial_t \chi_n  + \mu (\partial_t \chi_n)^2 
\quad\hbox{in } H, \hbox{ a.e. in } (0,T).
\label{22nbis}
\end{equation}

\section{A priori estimates}
\label{apriori}
\setcounter{equation}{0}

Now, we carry out the uniform estimates, independent of $n$, and useful to pass to the limit.

\noindent
\bf First a priori estimate.\quad\rm
Test \eqref{24n} by $\chi_n $ and \eqref{50} by $(1 -1/u_{n})$: this is formal but it can  be made rigorous 
by using a truncation of the function $r\mapsto (1-1/r)$ in a neighborhood of zero.  Then, we add the two
equalities by easily inferring
\begin{eqnarray*}
& & \frac\mu 2 \| \chi_n (t) \|_H^2 
+\frac1{2n} \|\nabla\chi_n (t)\|_H^2
+ \int_0^t\!\!\! \|\nabla\chi_n(s)\|_H^2 \,ds
+\int_0^t \langle \xi_n(s), \chi_n (s)\rangle \,ds\\
& & + \int_\Omega (u_n -\log u_n)(t) \,dx 
+ \int_0^t\!\!\!\int_\Omega (1+\chi_n) |\nabla\log u_n|^2 \,dxds
+ \gamma_n \int_0^t\!\!\!\int_\Gamma (1+\chi_n) (u_n -1) \,dxds\\
& & \leq \frac\mu 2 \| \chi_{0n} \|_H^2 
+\frac1{2n} \|\nabla\chi_{0n}\|_H^2
+\int_\Omega (u_{0n} -\log u_{0n}) \,dx\\
& & + \int_0^t\!\!\!\int_\Omega (h(\theta_n) - \log u_n)\chi_n \,dxds
- \int_0^t\!\!\!\int_\Omega u_n \nabla \chi_n \cdot \frac{\nabla u_n}{u_n^{\,2}} \,dxds.
\end{eqnarray*}
Owing to \eqref{25n}  and the monotonicity of  $\partial_{V,V'} J$, we have 
$\langle \xi_n, \chi_n\rangle \geq 0$ a.e. in $(0,T) $: indeed, note that 
the value $0$ corresponds to the minimum of $J$ (cf.~\eqref{defJ})
whence $0\in \partial_{V,V'} J (0)$. Moreover, since 
\begin{equation}
\label{posit}
0\leq \chi_{n}\leq 1 \ \hbox{ as well as } \ u_n >0 \ \hbox{ a.e. in $\, Q \, $ (and in $\,\Sigma\,$)}, 
\end{equation}  
and
\begin{equation}
r -\log r \geq \frac13 (r + |\log r| )\,, \quad r>0\,, 
\label{spirou}
\end{equation} 
it follows from \eqref{12} and \eqref{posit} that
$$
\int_0^t\!\!\!\int_\Omega (h(\theta_n) - \log u_n) \chi_n \,dxds \le c_h + \int_0^t\!\!\!\int_\Omega |\log u_n| \,dxds \le c_h + 3\ \int_0^t\!\!\!\int_\Omega (u_n - \log u_n) \,dxds ,
$$
while we remark that
\begin{eqnarray*}
- \int_0^t\!\!\!\int_\Omega u_n \nabla \chi_n \cdot \frac{\nabla u_n}{u_n^{\,2}} \,dxds
= - \int_0^t\!\!\!\int_\Omega \nabla \chi_n \cdot \nabla \log u_n \,dxds\\
\leq\frac12 \int_0^t\!\!\!\int_\Omega |\nabla\chi_n|^2 \,dxds + \frac12 \int_0^t\!\!\!\int_\Omega  |\nabla\log u_n|^2 \,dxds.
\end{eqnarray*}

We also recall \eqref{39}--\eqref{39bis} and finally apply 
the Gronwall lemma to deduce that
\begin{equation}
\begin{split}
\| \chi_n \|_{ L^\infty(0,T;H)\cap L^2(0,T;V)} & + \frac1{\sqrt n}\| \chi_n \|_{ L^\infty(0,T;V)} 
+ \| u_n \|_{ L^\infty(0,T;L^1(\Omega))} \\
& + \| \log u_n \|_{ L^\infty(0,T;L^1(\Omega))\cap L^2(0,T;V)}
+ \gamma_n \| u_n \|_{L^1(\Sigma)} \le C. 
\end{split}
\label{52}
\end{equation}

\noindent
\bf Second a priori estimate.\quad\rm
Test \eqref{24n} by $\partial_t \chi_n $: as the right hand side $h(\theta_n) - \log u_n$ of \eqref{24n} 
is bounded in $L^2(Q)$ thanks to \eqref{12} and \eqref{52}, by a standard computation we deduce that
\begin{equation}
\| \chi_n \|_{ H^1(0,T;H)\cap L^\infty(0,T;V)} + \frac1{\sqrt n}\| \partial_t\chi_n \|_{ L^2(0,T;V)} 
+ \| \hat\beta(\chi_n) \|_{ L^\infty(0,T;L^1(\Omega))} \leq C. 
\label{53}
\end{equation}

\noindent
\bf Third a priori estimate.\quad\rm
Rewrite \eqref{24n} as 
\begin{equation}
\label{var24n} 
\frac1n A (\partial_t \chi_n) + A \chi_n + \xi_n = f_n 
\end{equation} 
with 
$f_n:= h(\theta_n) - \log u_n - \mu\, \partial_t\chi_n $  already bounded in  $\ L^2(0,T;H)$.
 Then we can formally 
test \eqref{var24n}
by $A \chi_n$ (such estimate is performed rigorously in \cite[formulas~(173)--(175)]{chio}) 
and recover, with the help of well-known elliptic regularity results,
\begin{equation}
\| \chi_n \|_{L^2(0,T;W)} + \frac1{\sqrt n}\| A\chi_n \|_{ L^\infty(0,T;H)} \leq C
\label{54}
\end{equation}
provided \eqref{39bis} is used to control $n^{-1/2}\| A\chi_{0n} \|_{ H}.$ 

Now, thanks to \eqref{53}  and \eqref{54}, 
it turns out that  $\xi_n = f_n  - A (\partial_t \chi_n) / n - A \chi_n   $ fulfills
\begin{equation}
\| \xi_n \|_{L^2(0,T;V')} \leq C.
\label{54bis}
\end{equation}

\noindent
\bf Fourth a priori estimate.\quad\rm
We test \eqref{26n} by $p_n$ and arrive at 
\begin{eqnarray*}
&&\int_0^t\!\!\!\int_\Omega \frac{|p_n|^2}{1+\chi_n} \,dxds
+ \frac12 \|\nabla (1*p_n)(t)\|_H^2
+ \frac{\gamma_n}2 \int_\Gamma |(1*p_n)(t)|^2 \,dx\\
&&\hskip1cm\leq \|u_{0n}\|_{V'} \|(1*p_n)(t)\|_{V}
\leq \|u_{0n}\|_{V'} \int_0^t  \|p_n(s)\|_{H} \,ds 
+ \|u_{0n}\|_{V'} \|\nabla(1*p_n)(t)\|_{H}
\end{eqnarray*}
whence, as $(1+\chi_n)^{-1} \geq 1/2 $ a.e. in $Q$, using \eqref{39}, H\" older's and Young's inequalities 
we obtain
\begin{equation*}
\frac14\int_0^t\!\!\!\|p_n\|_H^2 \,ds
+ \frac14 \|\nabla (1*p_n)(t)\|_H^2
+ \frac{\gamma_n}2 \int_\Gamma |(1*p_n)(t)|^2 \,dx \leq C,
\end{equation*}
thus without exploiting the boundary contribution on the left hand side. 

Then, recalling also \eqref{27n}
we deduce that
\begin{equation}
\| u_n \|_{ L^2(0,T;H)} +\| 1*p_n \|_{ H^1(0,T;H)\cap L^\infty(0,T;V)} 
+ \sqrt{\gamma_n} \,\| 1*p_n \|_{ L^\infty(0,T;L^2(\Gamma))} \leq C. 
\label{55}
\end{equation}

\noindent
\bf Fifth a priori estimate.\quad\rm
We test \eqref{50} by $v(1+u_{n})^{-1}$ where $v$ is an arbitrary function in $W^{1, N+1}
(\Omega) \subset L^\infty (\Omega)$ and integrate only over $\Omega$. Note that, 
owing to \eqref{47} and the positivity of $u_n$, the  function $v(1+u_{n})^{-1}$ actually belongs to $V$.

Before proceeding, we also remind \eqref{52} and immediately deduce  that
\begin{equation}
\label{log1+un}
\| \log (1+u_n )\|_{ L^\infty(0,T;L^1(\Omega))\cap L^2(0,T;V)} \leq C
\end{equation}
from the  analogous estimate  on $\log u_n.$

Then, with the help of \eqref{50} and \eqref{posit} we obtain
\begin{eqnarray*}
&&\left|\int_\Omega \partial_t \log (1+u_n) \, v \,dx \right|\\
&&\hskip1cm = \left| \int_\Omega \nabla (u_n (1+\chi_n))\cdot \nabla \left( \frac{v}{1+u_n} \right) \,dx +
\gamma_n \int_\Gamma \frac{u_n}{1+u_n} (1+\chi_n) v \,dx
\right|\\
&&\hskip1cm \leq \left| \int_\Omega \left( (1+\chi_n)\nabla u_n + u_n \nabla \chi_n   \right) \cdot 
\left( \frac{1}{1+u_n} \nabla v  - \frac{v}{(1+u_n)^2} \nabla u_n \right) \,dx \right|
 + 2\gamma_n \int_\Gamma | v| \, dx\\
&&\hskip1cm \leq \int_\Omega \left( (1+\chi_n) \frac{|\nabla u_n|}{1+u_n} |\nabla v|
+ (1+\chi_n)|v| \frac{|\nabla u_n|^2}{(1+u_n)^2} \right) \,dx\\
&&\hskip1.5cm  {}+\int_\Omega \left(\frac{u_n}{1+u_n}  |\nabla \chi_n |\, |\nabla v|
+ \frac{u_n}{1+u_n} |v|\,  |\nabla \chi_n | \frac{|\nabla u_n|}{1+u_n} \right) \,dx + C \|v\|_{L^1(\Gamma)} .
\end{eqnarray*}
Consequently 
\begin{eqnarray*}
&&\left|\langle \partial_t \log (1+u_n) , v \rangle \right|\\
&&\hskip1cm \leq 2 \| \nabla\log (1+u_n )\|_{H} \| \nabla v \|_{H}
+2  \| \nabla\log (1+u_n )\|_{H}^2 \|  v \|_{L^\infty (\Omega)}\\
&&\hskip1.5cm  {}+ \| \nabla\chi_n \|_{H} \| \nabla v \|_{H}
+ \| \nabla\chi_n \|_{H}\| \nabla\log (1+u_n )\|_{H} \|  v \|_{L^\infty (\Omega)}
+ C \|v\|_{L^1(\Gamma)}\\
&&\hskip1cm \leq C \left( 1+\| \nabla\log (1+u_n )\|_{H}^2 + \| \nabla \chi_n \|_{H}^2 \right)
\|  v \|_{W^{1,N+1}(\Omega)}
\end{eqnarray*}
a.e. in $(0,T)$. Therefore, in view of \eqref{log1+un} and \eqref{53} we infer that 
\begin{equation}
\| \partial_t \log (1+u_n ) \|_{L^1(0,T;(W^{1,N+1}(\Omega))')}\leq C .
\label{dtlog1+un}
\end{equation}

\noindent
\bf Sixth a priori estimate.\quad\rm
We follow the argument devised by Boccardo and Gallou\" et in \cite{bg} 
and introduce the functions $\tau_k:\RR \to \RR$, $k\in\NN $,  defined by
$$
\tau_k(0)=0, \quad 
\tau'_k (r) =  \left\{
\begin{array}{ll}
\displaystyle
-1  \   &\hbox{if }  \ r \leq -(k+1)  
\\[0.1cm]
r+k \   &\hbox{if } \  -(k+1)\leq r  \leq -k  
\\[0.1cm]
0 \   &\hbox{if }  \ -k\leq r  \leq k  
\\[0.1cm]
r-k \   &\hbox{if } \  k \leq r  \leq k +1  
\\[0.1cm]
1  \   &\hbox{if } \  k+1\leq r 
\\
\end{array}
\right. .
$$
Note that $\tau_k\geq 0$, $|\tau'_k (r)|\leq 1,$ and $0\leq \tau''_k \leq 1$ a.e. in $\RR$. 
In view of \eqref{23n} and \eqref{2}, we have that
$$
\nabla e_n = (1+ \chi_n\theta_n h''(\theta_n))\nabla \theta_n  - (h(\theta_n)- \theta_n h'(\theta_n))\nabla \chi_n 
$$
whence 
$$
\displaystyle \nabla \theta_n = \frac{\nabla e_n + (h(\theta_n)- \theta_n h'(\theta_n))
\nabla \chi_n }{(1+ \chi_n\theta_n h''(\theta_n)) },
$$
the denominator being uniformly positive (in fact, bounded from below by $1/c_h$) due to \eqref{13}. 

Then, taking $\tau'_k(e_n)$ as test function in \eqref{22n}, and integrating in time, we have that
\begin{eqnarray}
&&\int_\Omega \tau_k(e_n(t)) \,dx
+ \int_0^t\!\!\!\int_\Omega \tau''_k (e_n)  \nabla e_n \cdot \nabla\theta_n \,dxds \nonumber\\
&&\hskip1cm \leq  \int_\Omega \tau_k(e_{0n}) \,dx + \int_0^t\!\!\!\int_\Omega \tau'_k (e_n) 
\left(\mu (\partial_t \chi_n)^2  - h(\theta_n)\partial_t \chi_n  \right) \,dxds .
\label{phi1}
\end{eqnarray}
With the help of Young's inequality  and \eqref{12}--\eqref{13}, we infer that
\begin{eqnarray*}
&&\int_0^t\!\!\!\int_\Omega \tau''_k (e_n)  \nabla e_n \cdot \nabla\theta_n \,dxds \\
&&\hskip1cm \geq
\int_0^t\!\!\!\int_\Omega \tau''_k (e_n)  \frac{|\nabla e_n|^2}{(1+ \chi_n\theta_n h''(\theta_n)) } \,dxds
+ \int_0^t\!\!\!\int_\Omega \tau''_k(e_n) \frac{(h(\theta_n)- \theta_n h'(\theta_n))}{(1+ \chi_n\theta_n h''(\theta_n)) } \nabla e_n \cdot \nabla \chi_n \,dxds\\
&&\hskip1cm\geq
\frac12\int_0^t\!\!\!\int_\Omega \tau''_k (e_n)  \frac{|\nabla e_n|^2}{(1+ \chi_n\theta_n h''(\theta_n)) } \,dxds
-\frac12 \int_0^t\!\!\!\int_\Omega \frac{(h(\theta_n)- \theta_n h'(\theta_n))^2}{(1+ \chi_n\theta_n h''(\theta_n)) } |\nabla \chi_n|^2 \,dxds\\
&&\hskip1cm\geq \frac1{2c_h}\int_0^t\!\!\!\int_\Omega \tau''_k (e_n)  |\nabla e_n|^2 \,dxds
- C \int_0^t\!\!\!\|\nabla \chi_n\|_H^2 \,ds.
\end{eqnarray*}
On the other hand, we observe that
$$
\int_0^t\!\!\!\int_\Omega \tau'_k (e_n) \left(\mu\, (\partial_t \chi_n)^2 - h(\theta_n)
\partial_t \chi_n  \right)dxds \leq C \int_0^t\!\!\! 
\left( 1 + \|\partial_t \chi_n\|_H^2 \right) ds.
$$
Hence, since $\tau_k(e_{0n})\leq |e_{0n}|$ and \eqref{39} holds, it follows from \eqref{phi1} and \eqref{53} that
\begin{equation}
\int_\Omega \tau_k(e_n(t)) \,dx + \int_0^t\!\!\!\int_\Omega \tau''_k (e_n)  
|\nabla e_n|^2 \,dxds \leq C . \label{phi2}
\end{equation}
Since $|r|\le 1+\tau_0(r)$ for $r\in\RR$, it follows in particular from \eqref{phi2} with $k=0$ that 
\begin{equation}
\|e_n \|_{L^\infty (0,T; L^1(\Omega))} \leq \sup_{t\in (0,T)} \int_\Omega (1+\tau_0(e_n(t))) \,dx 
\leq C. 
\label{phi3}
\end{equation}
Now, if we define $D_k:= \{(x,t)\in Q : \ k\leq |e_n(x,t) |\leq k+1\} $, \eqref{phi2} yields  
$$ 
\int\!\!\!\int_{D_k}  |\nabla e_n|^2 \,dxds \leq C \;\;\;\mbox{ for all }\;\;\; k\in \NN.
$$
Taking $1\leq q <2,$ we infer from the H\"older inequality that
\begin{eqnarray*}
\int\!\!\!\int_{Q}  |\nabla e_n|^q \,dxds & = & \sum_{k=0}^\infty \int\!\!\!\int_{D_k}  |\nabla e_n|^q \,dxds \\
& \leq & \sum_{k=0}^\infty \left( \int\!\!\!\int_{D_k}  |\nabla e_n|^2 \,dxds \right)^{\!\!q/2} |D_k|^{(2-q)/2} \leq C
\sum_{k=0}^\infty |D_k|^{(2-q)/2}.
\end{eqnarray*}
Moreover, let $q< (N+1)/(N+2)$ and $r = (N+1)q/N$. As $k \leq |e_n |$ a.e. in $D_k$, we have 
$$\displaystyle 
|D_k| \leq \frac1{k^r} \int\!\!\!\int_{D_k}  |e_n|^r \,dxds \quad \hbox{for }\ k\geq 1
$$
and, using \eqref{phi3} to control the integral on $D_0$, we obtain
\begin{eqnarray*}
\int\!\!\!\int_Q |\nabla e_n|^q \,dxds \leq 
C \left( 1 + \sum_{k=1}^\infty \frac1{k^{r(2-q)/2}} \left( \int\!\!\!\int_{D_k}  |e_n|^r \,dxds \right)^{\!\!(2-q)/2} \right)\\
\leq C \left( 1 + \left(\sum_{k=1}^\infty \frac1{k^{r(2-q)/q}} \right)^{\!\!q/2}  \left( \sum_{k=1}^\infty 
\int\!\!\!\int_{D_k}  |e_n|^r \,dxds \right)^{\!\!(2-q)/2} \right)\\
\leq C \left( 1 + \|e_n\|_{L^r(Q)}^{r(2-q)/2} \left(\sum_{k=1}^\infty \frac1{k^{r(2-q)/q}} \right)^{\!\!q/2}   \right).
\end{eqnarray*}
Let us point out that the series above converges if $\displaystyle \frac{r(2-q)}q >1 $, which follows 
from our choices of $r$ and $q$. Therefore, it results that 
\begin{equation}
\int\!\!\!\int_Q |\nabla e_n|^q \,dxds \leq C \left( 1 + \|e_n\|_{L^r(Q)}^{r(2-q)/2} \right).
\label{phi4}
\end{equation}
Now, by H\"older's inequality we deduce that
$$
\|e_n (t)\|_{L^r(\Omega)}^r \leq \|e_n (t)\|_{L^1(\Omega)}^\alpha \|e_n (t)\|_{L^{q^*}(\Omega)}^{q^*(1-\alpha)}
$$
with $ q^* = Nq/(N-q)$, $\alpha + (1-\alpha)q^*= r$, i.e., $1-\alpha= (N-q)/N = q/q^*$ if $N=2, 3$ and $ q^*= \infty$, $\alpha=1$, $q^*(1-\alpha)=r-1$ if $N=1$. We use the bound \eqref{phi3} to derive
$$
\|e_n (t)\|_{L^r(\Omega)}^r \leq C \|e_n (t)\|_{L^{q^*}(\Omega)}^{q}
$$
and, by Sobolev's embeddings, \eqref{phi3}, and \eqref{phi4}, we have 
\begin{eqnarray*}
 \|e_n\|_{L^r(Q)}^{r} \leq C \int_0^T  \|e_n (t)\|_{L^{q^*}(\Omega)}^{q} dt
 \leq C \left( 1+  \int_0^T  \|e_n (t) -\langle e_n(t),1 \rangle \|_{L^{q^*}(\Omega)}^{q} dt\right)\\
 \leq C \left( 1+  \int_0^T  \|\nabla e_n (t)   \|_{L^{q}(\Omega)}^{q} dt \right)
 \leq C \left( 1 + \|e_n\|_{L^r(Q)}^{r(2-q)/2} \right)
\end{eqnarray*}
from which a bound on $\|e_n\|_{L^r(Q)}$ follows since $2-q<2$. Recalling \eqref{phi4} and the fact that
$q$ varies in $[1,(N+2)/(N+1))$ and consequently $r$ in  $[1,(N+2)/N)$, we have thus established that 
\begin{equation}
\hbox{for all } \ 1\leq r <\frac{N+2}{N}, \  1\leq q < \frac{N+2}{N+1} \ \hbox{ there holds 
}\  \|e_n\|_{L^r(Q)\cap L^q(0,T;W^{1,q}(\Omega))} \leq C ,\quad 
\label{phi5}
\end{equation}
where the constant $C$ depends on $r$ and $q$ as well. 

In view of \eqref{12}, \eqref{13}, \eqref{23n}, and \eqref{52}, it is clear that $\theta_n$
satisfies the same estimates: 
\begin{equation}
\hbox{for all } \ 1\leq r <\frac{N+2}{N}, \  1\leq q < \frac{N+2}{N+1} \ \hbox{ there holds 
}\  \|\theta_n\|_{L^r(Q)\cap L^q(0,T;W^{1,q}(\Omega))} \leq C.
\label{spip}
\end{equation}
Thus, by comparison in \eqref{22n} and owing to the embedding of $L^1(\Omega)$ in $(W^{1,q'}(\Omega))'$, it turns out that
\begin{equation}
\| \partial_t e_n \|_{L^1(0,T;(W^{1,q'}(\Omega))')}\leq C \quad \hbox{ with } \ q'=\frac{q}{q-1} > N + 2 \quad \mbox{ for all }\;\;\; 1\leq q < \frac{N+2}{N+1}.
\label{phi6}
\end{equation}

This concludes the derivation of uniform estimates for our problem. 

\section{Passage to the limit as $n\to \infty$}
\label{limit}
\setcounter{equation}{0}

Taking advantage of the bounds in \eqref{posit}--\eqref{53}, \eqref{54}--\eqref{55} and \eqref{phi5}--\eqref{phi6},
we are in a position to pass to the limit as $n\to\infty$ by weak and weak star compactness arguments. Hence, possibly taking a subsequence (still labelled by $n$), it turns out that there are functions $e,\,\theta,\,\chi,\,\xi,\,u,\,p$ such that
\begin{eqnarray}
&e_n \to e \ \hbox{ and } \ \theta_n \to \theta \quad \hbox{ weakly in }
 L^r(Q)\cap L^q(0,T;W^{1,q}(\Omega)),&\nonumber \\
&\hskip4cm\hbox{for all $\,\displaystyle  1< r <\frac{N+2}{N} \,$ and $\,\displaystyle  1< q < \frac{N+2}{N+1}\, $,} \label{58}\\[0.2cm]
&\chi_n \to \chi \quad \hbox{ weakly star in } H^1(0,T;H)\cap L^\infty(0,T;V)\cap L^2(0,T;W)\cap L^\infty (Q),  &\label{60}\\
&\displaystyle \frac1n \partial_t \chi_n \to 0 \quad \hbox{ strongly in } L^2(0,T;V),  &\label{61}\\
&\displaystyle \xi_n \to \xi \quad \hbox{ weakly in } L^2(0,T;V'),  &\label{62}\\
&\displaystyle u_n \to u \quad \hbox{ weakly in } L^2(Q),  &\label{63}\\
&1*p_n \to 1*p \quad \hbox{ weakly star in } H^1(0,T;H)\cap L^\infty(0,T;V).  &\label{64}
\end{eqnarray}
We aim to prove that $(e,\theta,\chi,\xi,u,p)$ is a weak  solution to problem \eqref{1}--\eqref{8}
fulfilling the conditions stated in Theorem~\ref{existence}.

Thanks to \eqref{58} and to the boundedness \eqref{phi6} of $(\partial_t e_n ) $ in $L^1(0,T;(W^{1,q'}(\Omega))')$, it follows from
the compactness result in \cite[Sect.~8, Cor.~4]{Simon} that

\begin{equation}
e_n \to e  \quad \hbox{ strongly in } L^q(Q) \;\;\;\mbox{ for all }\;\;\; 1 \le q < \frac{N+2}{N+1}, \label{65}
\end{equation}
and this convergence actually takes place in $L^q(Q)$ for all $q\in [1,(N+2)/N)$ by interpolation with \eqref{58}. In addition, by the Ascoli theorem and the Aubin-Lions lemma~\cite[p.~58]{Lions} we also have
\begin{eqnarray}
&\chi_n \to \chi \quad \hbox{ strongly in } C^0([0,T];H)\cap L^2(0,T;V),  &\label{66}\\
&1*p_n \to 1*p \quad \hbox{ strongly in } C^0([0,T];H).  &\label{67}
\end{eqnarray}
Now, \eqref{65}--\eqref{66} and the Lipschitz continuity of $\psi^{-1} $ with respect to both 
variables (cf.~\eqref{psi2}) imply the strong convergence of $\theta_n = \psi^{-1}(e_n, \chi_n) $
to $\psi^{-1}(e, \chi)$ in $L^q(Q)$ for $q\in [1,(N+2)/N)$; hence, \eqref{58} and the uniqueness of the (weak) limit entail that
\begin{equation}
\theta_n = \psi^{-1}(e_n, \chi_n) \ \to \ \theta= \psi^{-1}(e, \chi) \quad \hbox{ strongly in } L^q(Q) \;\;\;\mbox{ for }\;\;\; 1 \le q < \frac{N+2}{N}, 
\label{68}
\end{equation}
and \eqref{23} follows. By virtue of \eqref{63}--\eqref{64}, \eqref{convgamma}, and \eqref{39}, passing to the limit in 
\eqref{26n} we recover \eqref{26} for $u$ and $1*p$. We point out that equality \eqref{27} is a consequence of  
\eqref{27n} since $(u_n (1+\chi_n))$ weakly converges to $u (1+\chi)$  in $L^1(Q)$ by \eqref{63} and \eqref{66}.

In view of \eqref{log1+un}--\eqref{dtlog1+un}, we use once more \cite[Sect.~8, Cor.~4]{Simon} to deduce that,
possibly taking a subsequence,
\begin{equation}
\log(1+u_n) \to \zeta \quad \hbox{ weakly in } \, L^2(0,T;V), \, \hbox{ strongly in } \, L^2(Q) \, \hbox{ and a.e. in }\, Q
\label{phi7}
\end{equation}
for some $\zeta \in L^2(0,T;V).$ Hence, it turns out that
$$ u_n = \exp (\log(1+u_n)) - 1 \ \to \ e^\zeta - 1 \quad \hbox{ a.e. in } Q ,$$
which, combined with \eqref{63}, yields $u= e^\zeta - 1 $ and 
\begin{equation}
u_n \to u \quad \hbox{ a.e. in } \, Q \, \hbox{ and strongly in } \, L^s(Q) \, \hbox{ for all } \ 1 \leq s <2
\label{phi8}
\end{equation}
(see, e.g., \cite[p.~12]{Lions}). Observe next that $\log u $ is well defined a.e. and integrable: indeed, 
by the Fatou lemma one sees that
$$ 
\int\!\!\!\int_Q |\log u| \, dxds = \int\!\!\!\int_Q \liminf_{n\to\infty}
 |\log u_n| \, dxds
\leq 
\liminf_{n\to\infty} \int\!\!\!\int_Q |\log u_n| \, dxds   \leq C.
$$
Then, \eqref{phi8} and \eqref{52} entail 
\begin{equation}
\begin{split}
\log u_n \ \to \ \log u &\hbox{ weakly in } \, L^{8/3}(Q)\cap L^2(0,T;V) \\
&\hbox{ and strongly in } \, L^s(Q) \, \hbox{ for all } \ 1 \leq s <8/3,
\end{split}
\label{phi9}
\end{equation}
where $8/3$ is chosen in  a way that $L^\infty (0,T;L^1 (\Omega)) \cap L^2(0,T;V) \subset L^{8/3} (Q)$ 
for $N\in\{1,2,3\}$. 
We also recall \eqref{68} and the Lipschitz continuity and boundedness of $h$, which imply
\begin{equation}
h(\theta_n )\ \to \ h(\theta) \quad\hbox{ strongly in } \, L^s(Q) \, \hbox{ for all } \ 1 \leq s < \infty .
\label{68bis}
\end{equation}
Then, using \eqref{60}--\eqref{62} we can pass to the limit in \eqref{24n} and obtain \eqref{24}. From \eqref{25n} and the inequality
$$
\limsup_{n\to\infty} \int_0^T \langle \xi_n (t), \chi_n(t)\rangle dt \leq
\int_0^T \langle \xi (t), \chi (t)\rangle dt,
$$
which is a consequence of \eqref{62} and \eqref{66}, the inclusion \eqref{25} results from 
\cite[Lemma~1.3, p.~42]{Barbu} for $\partial_{V,V'} J$ induces a maximal monotone operator also from $L^2(0,T;V)$ 
to $L^2(0,T;V')$. Note in particular that  $\chi $ satisfies \eqref{31bis} (look at \eqref{60} and use interpolation for the continuity  from $[0,T] $ to $V$) and $\xi$ is in $L^2(0,T;H)$, that is \eqref{32bis},
due to  a comparison of terms in \eqref{24} (cf. also Remark~\ref{rem-dopo-def}).

Now, it is our intention to show that
\begin{equation}
\partial_t\chi_n \to \partial_t \chi  \quad \hbox{ strongly in } \, L^2(Q) ,
\label{72}
\end{equation}
in order to be able to pass to the limit in \eqref{22n}. 
Test \eqref{24n} by $\partial_t\chi_n$ and observe that
\begin{eqnarray}
&&\mu\, \limsup_{n\to\infty} \int\!\!\!\int_Q |\partial_t \chi_n|^2 \,dxds\nonumber\\
&&\hskip1cm = \limsup_{n\to \infty} \left( - \frac1n \int\!\!\!\int_Q |\nabla (\partial_t \chi_n)|^2 \,dxds
- \frac12 \|\nabla \chi_n(T)\|_H^2  - J(\chi_n(T)) \right.\nonumber\\
&&\hskip3.5cm \left. {}+  \frac12 \|\nabla \chi_{0n}\|_H^2 
{}+ J(\chi_{0n}) + \int\!\!\!\int_Q (h(\theta_n ) - \log u_n) \partial_t \chi_n  \,dxds \right)\nonumber\\
&&\hskip1cm \leq - \liminf_{n\to \infty} \left( 
 \frac12 \|\nabla \chi_n(T)\|_H^2 + J(\chi_n(T)) \right.\nonumber\\
&&\hskip3.5cm \left. {}-  \frac12 \|\nabla \chi_{0n}\|_H^2 
{}-  J(\chi_{0n}) {}- \int\!\!\!\int_Q (h(\theta_n ) - \log u_n) \partial_t \chi_n \,dxds \right)\nonumber\\
&&\hskip1cm \leq - \left( 
 \frac12 \|\nabla \chi (T)\|_H^2 + J(\chi(T)) 
-  \frac12 \|\nabla \chi_{0}\|_H^2 
-  J(\chi_{0}) - \int\!\!\!\int_Q (h(\theta ) - \log u) \partial_t \chi \,dxds \right)\hskip1.5cm \label{semic}
\end{eqnarray}
by the weak lower semicontinuity of the semi-norm in $V$ and the convex function $J,$ as $\chi_n(T) \to \chi(T)$ weakly in V due to \eqref{60} and \eqref{66}, and \eqref{39}, \eqref{39ter}, \eqref{phi9} and \eqref{68bis} as well. Now, we can exploit \eqref{24} and conclude that the last line of \eqref{semic} is nothing but  $\|\partial_t \chi \|^2_{L^2(Q)}$, as one can see by testing \eqref{24} by $\partial_t \chi $ and taking advantage of \eqref{25} and \cite[Lemma~3.3, p.73]{Brezis}. Hence, we have that 
$$
\limsup_{n\to\infty} \|\partial_t \chi_n\|^2_{L^2(Q)} \leq \|\partial_t \chi \|^2_{L^2(Q)}
$$
which, along with the weak convergence \eqref{60} of $(\partial_t \chi_n)$ to  $\partial_t \chi $ in $L^2(Q)$, 
ensures \eqref{72}. 

At this point, we can take the limit in \eqref{22n}. Thanks to \eqref{39},  \eqref{58}, \eqref{68}, \eqref{68bis}--\eqref{72} it turns~out~that
$$
e_n (t) = e_{0n} - \int_0^t \left( A \theta_n  + h(\theta_n) \partial_t  
\chi_n -  \mu (\partial_t \chi_n)^2 \right)(s)\, ds
$$
converges in a suitable way to
$$
e (t) = e_{0} - \int_0^t \left( A \theta  + h(\theta) \partial_t \chi  
-  \mu\, (\partial_t \chi )^{\,2 }\right) (s)\, ds ,
$$ 
whence the initial condition \eqref{28} for $e$ and \eqref{22}. Note in particular that $e$ is continuous from  $[0,T]$ to $(W^{1,q'}(\Omega))' $ and satisfies \eqref{30} and \eqref {29bis} for all $r,\,q $ fulfilling $1\leq r <(N+2)/N$,  $1\leq q < (N+2)/(N+1)$. In fact, the latter comes as a consequence of the bounds \eqref{phi5} coupled with the convergences \eqref{58}. About the 
$L^\infty(0,T;L^1(\Omega))$-regularity of $e$, \eqref{65} implies that $ e_n \to e $ in $L^1 (0,T;L^1(\Omega))$, whence
$$
\| e_n(t)\|_{L^1(\Omega)} \ \to \ \| e (t)\|_{L^1(\Omega)} \quad\hbox{ for a.e. } \, t\in(0,T), \, \hbox{ at least for a subsequence.}
$$
Then, recalling \eqref{phi3} we infer that
$$ \| e (t)\|_{L^1(\Omega)}\leq \sup_{n\to \infty}  \| e_n\|_{L^\infty(0,T;L^1(\Omega))} \leq C,$$
and consequently $e \in L^\infty(0,T;L^1(\Omega)).$ The same property can be deduced for $u$ and $\log u$
(cf. \eqref{phi8} and \eqref{phi9}) so that \eqref{33bis} holds as well. The proof of Theorem~\ref{existence} is then complete.

\section{Complementary results}
\label{complement}
\setcounter{equation}{0}

This section is devoted to the proofs of the other results stated in Section~\ref{mainres}. 
We start with Theorem~\ref{posit-temp} ensuring the positivity of the variable $\theta$.

\noindent
\bf Proof of Theorem \ref{posit-temp}.\quad\rm
Let us go back to our approximating problem \eqref{22n}--\eqref{28n}, recalling  in particular Proposition~\ref{app-sol}. By \cite[Theorem~2.6]{chio} we have 
\begin{equation}
\log \theta_n  \in L^\infty(0,T;L^1(\Omega))\cap L^2(0,T;V), \quad 
\frac{\partial_t\chi_n}{\sqrt{\theta_n}} 
\in L^2(Q) 
\label{reg1n}
\end{equation}
whenever $\log \theta_{0n} \in L^1(\Omega)$, besides the other requirements \eqref{37}--\eqref{38} 
on the initial data. Then, here we have not only to guarantee that $\log \theta_{0n} \in L^1(\Omega)$ 
for all $n\in\NN$, but we also need the uniform boundedness property
\begin{equation}
\|\log\theta_{0n}\|_{L^1(\Omega)} \leq C
\label{theta0-unif}
\end{equation}
along with \eqref{39}--\eqref{39bis} to prove \eqref{reg1} for the limiting functions. 

By \eqref{psi1}, \eqref{35}, and \eqref{hyptheta0} we know that both $\theta_0$ and $\log \theta_0$ are in 
$L^1(\Omega)$. 
Then we could, for instance, introduce $\eta_0:= \sqrt{\theta_0}$, observe that $\eta_0\in H$ and
$\displaystyle \log \eta_0 = \frac12 \log \theta_0 \in L^1(\Omega)$, and consequently approximate 
$\eta_0$ by taking the solution $\eta_{0n}$ to (cf.~\eqref{uzeron})
\begin{equation}
\label{etazeron}
\eta_{0n} + \frac1n A \eta_{0n} = \eta_0 + \frac{1}n \quad \hbox{in } H,
\end{equation}
then setting $\theta_{0n} := \eta_{0n}^{\, 2}$. By this choice, 
it  turns out that $\eta_{0n} \in W $ (whence $\theta_{0n}\in V$) and
$\eta_{0n} \geq 1/n $ a.e. in $\Omega$. Moreover, $\eta_{0n} \to \eta_0 $ in $H $
so that $\theta_{0n} \to \theta_0 $ in $L^1(\Omega) $ as $n\to\infty$. Then the related 
sequence $e_{0n} := \psi(\theta_{0n}, \chi_{0n}) $ satisfies $\eqref{39}_1$, i.e., converges 
to $e_0$ in $ L^1(\Omega)$ thanks to the properties of $\psi$ (see \eqref{psi1}--\eqref{psi2}).
Furthermore, if we test \eqref{etazeron} by $(1 -1/\eta_{0n})$, use \eqref{spirou}, and proceed as in the computation 
following \eqref{uzeron}, we arrive at
\begin{eqnarray*}
&&\frac16 \int_\Omega  |\log \theta_{0n}| \,dx \leq \frac13 \int_\Omega (|\eta_{0n}| + |\log \eta_{0n}|) \, dx\\
&&\hskip1cm {}\leq \int_\Omega (\eta_{0n} - \log \eta_{0n}) \, dx\leq  \int_\Omega (\eta_{0} - \log \eta_{0}) \,dx + |\Omega|,
\end{eqnarray*}
whence \eqref{theta0-unif} is satisfied.

Now, we can deal with the approximating problem as in Proposition~\ref{app-sol} 
and, in particular, consider equation \eqref{22nbis} which is an alternative formulation of \eqref{22n}. We formally
test \eqref{22nbis} by  $-1/\theta_n$ and integrate by parts: we easily infer that 
\begin{eqnarray*}
&&\int_\Omega (- \log\theta_n) (t) \,dx +\int_0^t\!\!\!\int_\Omega \frac{|\nabla\theta_n|^2}{\theta_n^{\,2}} \,dxds
+\mu \int_0^t\!\!\!\int_\Omega \frac{|\partial_t\chi_n|^2}{\theta_n} \,dxds \\
&&\hskip1cm {}\leq \int_\Omega ( - \log \theta_{0n}) \,dx + \int_\Omega (h'(\theta_n) \chi_n) (t) \,dx -
\int_\Omega h'(\theta_{0n}) \chi_{0n} \,dx .
\end{eqnarray*}
Then we add the term $\displaystyle \int_\Omega  \theta_n (t) \,dx$ to both sides of the above inequality 
and recall \eqref{12}, \eqref{posit}, \eqref{spirou}, and \eqref{spip} to deduce that
\begin{eqnarray*}
&&\frac13 \int_\Omega (|\theta_{n}(t)| +  |\log\theta_{n} (t)|) \,dx +\int_0^t\!\!\!\int_\Omega |\nabla\log \theta_n|^2 \,dxds
+\mu \int_0^t\!\!\!\int_\Omega \frac{|\partial_t\chi_n|^2}{\theta_n} \,dxds\\
&&\hskip1cm {}\leq \int_\Omega ( |\log \theta_{0n}| + |h'(\theta_{0n}) \chi_{0n}|) \,dx
+ \int_\Omega ( \theta_n +  h'(\theta_n) \chi_n) (t) \,dx \leq C,
\end{eqnarray*}
whence we find out the estimate 
\begin{equation}
\|\log \theta_n \|_ {L^\infty(0,T;L^1(\Omega))\cap L^2(0,T;V)} 
+  \left\|\frac{\partial_t\chi_n}{\sqrt{\theta_n}}\right\|_ {L^2(Q)} \leq C .
\end{equation}

Next, we can pass to the limit as  $n\to\infty$. In view of \eqref{68} and \eqref{72}, 
we have that  $\theta_n \to \theta$ and $ \partial_t\chi_n \to \partial_t \chi  $ a.e. in $Q$, at least 
for a subsequence. Thus, similarly as in \eqref{phi9} we infer that $\theta>0$ a.e. in $Q$ and
\begin{eqnarray*}
\log \theta_n \ \to \ \log \theta &&\hbox{ weakly in } \, L^{8/3}(Q)\cap L^2(0,T;V)\nonumber\\
&&\qquad\hbox{ and strongly in } \, L^s(Q) \, \hbox{ for all } \ 1 \leq s <8/3, \\
\frac{\partial_t\chi_n}{\sqrt{\theta_n}} \ \to \ \frac{\partial_t \chi }{\sqrt{\theta}} &&\hbox{ weakly in } \, L^2(Q)
\, \hbox{ and strongly in } \, L^s(Q) \, \hbox{ for all } \ 1 \leq s <2.
\end{eqnarray*}
At this point, it is a standard matter to recover \eqref{reg1} and Theorem~\ref{posit-temp} 
is completely proven.\QED

The next step is showing the stability and decay estimate in the case $\gamma>0$.

\noindent
\bf Proof of Theorem \ref{decay}.\quad\rm
We first note that for a fixed $\chi $ fulfilling \eqref{31}, \eqref{25bis} and \eqref{31bis} (as it is for 
our  solution), there exists  a unique pair $(u,p)$ solving \eqref{26}--\eqref{27} and such that
\begin{eqnarray}
& u\in C^0([0,T];V') \cap L^2(0,T;H), & \label{33ter}\\
& p \in  L^2(0,T;H)   \ \hbox{ with } \  1*p \in C^0([0,T];V)  & \label{34ter}
\end{eqnarray}
for all $T>0.$ Now, we proceed formally by letting $u$ and $p$ be smooth enough in order to perform our computations:
for instance, a rigorous proof can be obtained by suitably regularizing $u_0$ and the fixed $\chi$ and then passing to the limit. Hence, we consider the solution $\zeta \in C^0([0,+\infty); V)$ to 
\begin{equation} 
A_\gamma \zeta(t) = u (t) \quad\hbox{in } V', \hbox{ for all } t\geq 0  \label{zeta}
\end{equation}
and perform a formal estimate by testing \eqref{zeta} by $\partial_t \zeta (t)$. Thus, with the help
of the time differentiated version of \eqref{26}, 
$$
\partial_t u (t) + A_\gamma p (t) = 0 , \quad t \in (0,T), 
$$ 
we easily obtain 
\begin{eqnarray}
\frac12 \frac{d}{dt} \langle  A_\gamma \zeta(t), \zeta (t) \rangle = 
\langle  A_\gamma \zeta_t (t) , \zeta (t) \rangle = 
\langle  \partial_t u (t), \zeta (t) \rangle = - \langle  A_\gamma p (t), \zeta (t) \rangle \qquad \nonumber\\
{}= - \langle  A_\gamma \zeta (t), p (t) \rangle =  
- \langle  u (t), p (t) \rangle  . \label{iden} 
\end{eqnarray}
Owing to \eqref{27} and setting  
$$
 \Phi_1(t) := \frac12 \|\nabla \zeta (t)\|_H^2 + \frac{\gamma}2  \int_\Gamma |\zeta (t)|^2 \,dx, \quad 
\Phi_2 (t) := \int_0^t\!\!\!\int_\Omega (1+\chi) u^2 \,dxds, 
$$
from \eqref{iden} it turns out that the continuous and nonnegative function $\Phi_1 + \Phi_2  $ is constant, or 
equivalently, 
\begin{equation}
 \frac{d\Phi_1}{dt}(t)  +  \int_\Omega \left( (1+\chi) u^2\right)(t) \,dx =0 \quad \hbox{ for a.e. } \, t \in  (0,+\infty). 
\label{null}
\end{equation}
Now, since \eqref{zeta} implies 
$$
\|\nabla \zeta (t)\|_H^2 + \gamma  \int_\Gamma |\zeta (t)|^2 \,dx\leq 
 \|u (t)\|_{V'}  \|\zeta (t)\|_V  \leq C\, \|u (t)\|_{V'} \left(\|\nabla \zeta (t)\|_H^2 
+ \gamma  \int_\Gamma |\zeta (t)|^2 \, dx \right)^{1/2}
$$
and, recalling that
$$
\|u (t)\|_{V'}^2 \leq C \|u (t)\|_{H}^2 \leq 
C \int_\Omega \left( (1+\chi) u^2\right)(t) \,dx \, 
$$ 
for a.e. $t \in  (0,+\infty)$, it is not difficult to see that \eqref{null} entails
\begin{equation}
 \frac{d}{dt} \Phi_1(t)  +  \alpha  \Phi_1(t) \leq 0
\label{null2}
\end{equation}
for some positive constant $\alpha.$ Then the Gronwall lemma yields  $\Phi_1(t) \leq \Phi_1(0) e^{-\alpha t},$
which means exponential decay for both $\|\zeta (t)\|_V $ and $\|u (t)\|_{V'}$ since (cf.~\eqref{zeta}) 
$$ 
\|u (t)\|_{V'} = \|A_\gamma \zeta (t)\|_{V'} \leq C \|\zeta (t)\|_V.
$$
Moreover, a comparison in \eqref{26} shows that $(1*p) (t)$ converges to $A_\gamma^{-1}(u_0)$ in $V$ (still at 
exponential rate) as $t\to +\infty$. Now, by integrating \eqref{null} with respect to time and using the nonnegativity of $\chi$, we have
$$  \| u\|_{L^2(0,t; L^2(\Omega))} = \int_0^t\!\!\!\int_\Omega  u^2 \,dxds\leq \Phi_2(t) \leq \Phi_1(0) \leq C \quad \mbox{ for all }\;\;\ t \in  (0,+\infty),$$
which completes the proof of Theorem~\ref{decay}.\QED

Finally, we prove the result on steady state solutions for the case $\gamma = 0$.

\noindent
\bf Proof of Theorem \ref{stationary}.\quad\rm
Here, we are looking for steady state solutions of problem~\eqref{1}--\eqref{8} with $\gamma = 0$: let us 
remind Definition~\ref{solution} and look for sextuples of functions $(e,\,\theta,\,\chi,\,\xi,\,u,\,p)$ which do not depend on time.

In view of \eqref{22}, we observe that $A\theta =0$ implies that $\theta$ is a constant function; the same can be concluded for $p$ from \eqref{26}. We now look for stationary $\chi$ solving \eqref{24}, which can be conveniently rewritten as 
\begin{equation}  
 A \chi + \omega = h(\theta) - \log p
\label{eq-chi}
\end{equation}
where (cf.~\eqref{4}) $\omega : = \xi - \log (1+\chi)$ satisfies
\begin{equation}
\omega \in \partial I_{[0,1]} (\chi)  \quad \hbox{ a.e. in }\, 
\Omega.
\label{eq-ome}
\end{equation} 
Since the right hand side of \eqref{eq-chi} is constant, two solutions 
$\chi_1, \, \chi_2 $ of \eqref{eq-chi} may differ just by a constant: in fact, if we test 
the difference of equations by $\chi_1 - \chi_2$, due to the monotonicity of $ \partial 
I_{[0,1]}$ we have that $ |\nabla (\chi_1 - \chi_2)| = 0 $ a.e. in $\Omega$. Therefore, we 
can observe that (recall that $p=u(1+\chi)$)
\begin{enumerate}
\item[(i)] if $ h(\theta) - \log p >0$ then $\chi = 1 $ and $\omega = h(\theta) - \log p$ 
solve \eqref{eq-chi}--\eqref{eq-ome}. According to the above discussion, any other solution $\chi $ should be a 
constant, but no constant different from $1$ is allowed by \eqref{eq-ome}. Therefore,  $\chi = 1 $ and $\omega = h(\theta) - \log p$ is the only solution to \eqref{eq-chi}--\eqref{eq-ome} in that case.
\item[(ii)]  if $ h(\theta) - \log p < 0$ then $\chi = 0 $ and $\omega = h(\theta) - \log p$ 
solve \eqref{eq-chi}--\eqref{eq-ome}. Any other solution $\chi $ should be a constant,
but no constant different from $0$ is allowed by \eqref{eq-ome}, so that $\chi = 0 $ and $\omega = h(\theta) - \log p$ is the only solution to \eqref{eq-chi}--\eqref{eq-ome} in that case.
\item[(iii)] if $ h(\theta) - \log p = 0$ then any constant  $\chi \in [0,1]$ solves 
\eqref{eq-chi}--\eqref{eq-ome} with $\omega = 0 $. 
\end{enumerate}
Then $\chi$ turns out to be a constant in all cases and consequently also $ e = \psi(\theta,
\chi)$ given by \eqref{23} is a constant. This concludes the proof of Theorem~\ref{stationary}.\QED


\section*{Acknowledgments}

The authors gratefully acknowledge financial support and kind hospitality of the 
Institut de Math\'e\-matiques de Toulouse, Universit\'e Paul Sabatier, and the
IMATI of CNR in Pavia. Credits also go to the ERC grant BioSMA (2008-2013) 
and the MIUR-PRIN Grant 2008ZKHAHN ``Phase transitions, hysteresis 
and multiscaling''.



\begin{thebibliography}{20}


\bibitem{Barbu}
V. Barbu,
``Nonlinear semigroups and differential equations in Banach spaces'',
Noord\-hoff,
Leyden,
1976.

\bibitem{bg}
L. Boccardo, T. Gallou\" et, 
Nonlinear elliptic and parabolic equations involving measure data,
{\it J. Funct. Anal.} {87} (1989) 149 - 169.

\bibitem{B1}
E. Bonetti, 
Global solution to a nonlinear phase transition model with dissipation,
{\it  Adv. Math. Sci. Appl.} {12} (2002) 355 - 376. 

\bibitem{bfl}
E. Bonetti, M. Fr\' emond, and C. Lexcellent, 
Hydrogen storage: modeling and analytical results, 
{\it Appl. Math. Optim.} {55} (2007) 31 - 59. 

\bibitem{BFF}
E. Bonetti, M. Fabrizio, and M. Fr\' emond, 
A first order phase transition with non-constant density, 
{\it J. Math. Anal.  Appl.} to appear (2011) doi: 10.1016/j.jmaa.2011.06.017. 

\bibitem{BFL}
G. Bonfanti, M. Fr\'emond, and F. Luterotti,
Global solution to a nonlinear system for irreversible phase changes,
{\it  Adv. Math. Sci. Appl.} {10} (2000) 1 - 24.

\bibitem{Brezis}
H. Brezis,
``Op\'erateurs maximaux monotones et semi-groupes de contractions
dans les espaces de Hilbert'',
North-Holland Math. Stud.
{\bf 5},
North-Holland,
Amsterdam,
1973.

\bibitem{chio}
E. Chiodaroli,
A dissipative model for hydrogen storage: existence and regularity results,
{\it Math. Methods Appl. Sci.} {34} (2011) 642 - 669.

\bibitem{CLSS}
P. Colli, F. Luterotti, G. Schimperna, and U. Stefanelli,
Global existence for a class of generalized systems for irreversible phase changes,
{\it NoDEA Nonlinear Differential Equations Appl.} {9} (2002) 255 - 276.

\bibitem{FPR}
E. Feireisl, H. Petzeltova, and E. Rocca,
Existence of solutions to a phase transition model with microscopic movements,
{\it Math. Methods Appl. Sci.} {32} (2009) 1345 - 1369.

\bibitem{Fremond}
M. Fr\'emond,
``Non-smooth Thermomechanics'',
Springer-Verlag, Berlin, 2002.

\bibitem{RF}
M. Fr\'emond, E. Rocca,
Well-posedness of a phase transition model with the possibility of voids,
{\it Math. Models Methods Appl. Sci.} {16} (2006) 559 - 586.

\bibitem{LSS}
Ph. Lauren\c cot, G. Schimperna, and U. Stefanelli,
Global existence of a strong solution to the one-dimensional full model for irreversible phase transitions,
{\it  J. Math. Anal. Appl.} (271) (2002) 426 - 442.

\bibitem{Lions}
J.L. Lions,
``Quelques m\'ethodes de r\'esolution des probl\`emes aux limites non
lin\'eaires'',
Dunod Gauthier--Villars,
Paris,
1969.

\bibitem{LS}
F. Luterotti, U. Stefanelli,
Existence result for the one-dimensional full model of phase transitions,
{\it  Z. Anal. Anwendungen} {21} (2002) 335 - 350.

\bibitem{Simon}
J. Simon,
{Compact sets in the space $L^p(0,T; B)$},
{\it Ann. Mat. Pura Appl.} {\bf 146} (1987) 65 - 96.

\end{thebibliography}
\end{document}